\documentclass[11pt,a4paper,oneside]{amsart}
\usepackage{amsmath}
\usepackage{amssymb}
\usepackage{latexsym}

\newcommand{\co}{\colon\thinspace}    %  Colon with correct spacing for maps.
\newcommand{\fnote}[1]{\footnote{\small sharp1}}
\newcommand{\inv}{^{-1}}              %  inverse
\newcommand{\Pom}{\mbox{${\mathcal A}_{\omega}$}}
\newcommand{\Po}{\mbox{${\mathcal A}_{0}$}}
\newcommand{\Pu}{\mbox{${\mathcal A}_{(u^+,u^-)}$}}

\newcommand{\To}{\mbox{${\mathcal T}_{0}$}}
\newcommand{\Ph}{\mbox{${\mathcal A}_{h}$}}
\newcommand{\Mh}{\mbox{${\mathcal M}_{h}$}}
\newcommand{\Mom}{\mbox{${\mathcal M}_{\omega}$}}

\newcommand{\vect}{\mbox{Vect}}

\newcommand{\supp}{\mbox{supp}}
\newcommand{\inter}{\mbox{Int}}

\def \R{\mathbb{R}}%        corps des reels
\def \Z{\mathbb{Z}}
\def \N{\mathbb{N}}

\newcommand{\al}{\alpha}
\newcommand{\om}{\omega}
\newcommand{\ga}{\gamma}
\newcommand{\ep}{\epsilon}
\newcommand{\gadot}{\dot{\gamma}}
\newcommand{\spt}{\mbox{supp}}

\newtheorem{theorem}{Theorem}
\newtheorem{proposition}[theorem]{Proposition}
\newtheorem{corollary}[theorem]{Corollary}
\newtheorem{definition[theorem]}{Definition}
\newtheorem{lemma}[theorem]{Lemma}
 
\newtheorem{conjecture}[theorem]{Conjecture} 

\title{On Aubry Sets and Mather's Action Functional}
\author{Daniel Massart}
\date{\today}
\begin{document}
\begin{abstract}
\noindent
We study Lagrangian systems on a closed manifold $M$.
We link the differentiability of Mather's $\beta$-function with the 
topological complexity of the complement of the Aubry set. 
As a consequence, when $M$ is a closed, orientable surface,
the differentiability of the $\beta$-function at a given homology class
is forced by the irrationality of the homology class. 
This allows us to prove the two-dimensional case 
of a conjecture by Ma\~{n}\'e. 
\end{abstract}

\maketitle

\section{Introduction}
We start by recalling some facts about Aubry-Mather theory.
Let $M$ be a smooth, closed, connected $n$-dimensional
manifold and $L$ be a Lagrangian on
the tangent bundle $TM$, that is, a $C^r, r\geq 2$ 
function on $TM$ which is convex
and superlinear when restricted to any fiber. 
The Euler-Lagrange equation then
defines a  flow $\Phi_t$ on $TM$, complete in the autonomous case.
Throughout this paper we assume 
$M$ to be endowed with a fixed Riemann metric, with respect to which we 
evaluate distances and norms in the tangent bundle ; our results do not 
depend on the  metric.
Denote by $\pi$ the canonical projection $TM \rightarrow M$. 

For $x,y \in M$ define $h_t(x,y)$ as the minimum, over all 
absolutely continuous curves $\ga \co [0,t] \rightarrow M$
with $\ga(0)=x,\, \ga(t)=y$, of $\int_0^t L(\ga,\gadot)ds$.
Then, by Fathi's weak KAM  theorem (\cite{Fathi1}) there exists
$c(L) \in \R$ such that 
$\liminf_{t\rightarrow \infty}(h_t(x,y) +c(L)t)$
is finite for every $x,y$.
This $\liminf$, originally defined in \cite{Mather3}, 
is called the Peierls barrier and denoted $h(x,y)$  
and $c(L)$ is Ma\~n\'e's critical value (see \cite{Mane3}). 
The Aubry set $\Po$ is then defined in \cite{Fathi2}
as the zero locus of $h$ restricted to the diagonal in $M \times M$. 
The canonical projection $\pi$ is a bi-Lipschitz homeomorphism between
$\Po $  and the set $\mbox{$\tilde{{\mathcal A}}_0$}$  
of velocity  vectors of  orbits in $\Po$ (Graph Property).
Furthermore $\mbox{$\tilde{{\mathcal A}}_0$}$ 
is compact and $\Phi_t$-invariant. 
 
Fathi's weak KAM  theorem  asserts that there exists
a Lipschitz function $u_+$ (resp. $u_-$)
such that $u_{\pm}(\ga(t))-u_{\pm}(\ga(0)) \leq  
\int_0^t (L +c(L))(\ga,\gadot)ds$ 
for every absolutely continuous path $\ga \co [0,t] \rightarrow M$,
which is written $u \prec L +c(L)$ for short, 
and such that for every $x \in M,\, t \in \R$ there is a
$C^1$ path $\ga \co [0,t] \rightarrow M$ with $\ga(0) = x$ 
(resp. $\ga(t) = x$) achieving equality. 
Such functions come in 
pairs, called conjugate pairs $(u_+, u_-)$ such that $u_+ \leq u_-$ with 
equality on $\Po$. 
Theorem 6 of \cite{Fathi2} asserts that
$h(x,y) = \sup \{u_-(y)-u_+(x)\}$, where the supremum is taken over 
conjugate pairs of weak KAM solutions.

For every closed 1-differential $\omega$, $L-\omega$ is a convex and 
superlinear Lagrangian, we sometimes
denote $\Pom$ its Aubry set $\Po (L-\omega)$.
Mather's $\alpha$-function is defined in \cite{Mather1} as
\[
 \alpha(\omega) = -\min \{ \int_{TM} (L-\omega)d\mu \co \mu \in {\mathcal M}
\} 
\]
where $\mathcal{M}$ is the set of closed 
measures  on $TM$, that is (see \cite{Bangert3})
the compactly supported probability measures $\mu$
on $TM$ such that $\int df \, d\mu =0$ for every $C^1$ function $f$ on $M$.
In other words, those are the measures with a well-defined homology class.
The measures achieving the minimum are invariant by the Euler-Lagrange 
flow $\Phi_t$ of $L$ (see \cite{Bangert3}). 
The quantity $\alpha$ defines a convex and superlinear 
function on $H^1(M, \R)$, twice the squareroot of which is also called stable 
norm when $L$ is a metric (see \cite{Massart} and the references therein). 
It is convex and superlinear
and its Fenchel transform is Mather's $\beta$-function on $H_1(M, \R)$, which 
is defined, for every real homology class $h$, as
\[
\beta(h) = \min \{ \int_{TM} (L)d\mu \co  \mu \in {\mathcal M},\, [\mu]=h
\}. 
\]
Let $\tilde{{\mathcal M}}_{\om}$ be  the closure in $TM$ 
of the union of the supports of  
measures in ${\mathcal M}$ achieving the minimum in the expression of $\al$.
Such measures are called 
$\om$-minimising measures, or just minimising measures if $[\om]=0$. 
We call Mather set of $L$ and $\omega$, 
and denote ${\mathcal M}_{\om}$ the projection 
$\pi(\tilde{\mathcal M}_{\om})$ ;
it is contained in $\Pom$ (\cite{Fathi1}).
In particular we call Mather set of $L$ the Mather set ${\mathcal M}_0$ 
corresponding to the zero cohomology class. 

For every $[\omega] \in H^1(M, \R)$ we call $F_{\omega}$ the 
maximal face of the epigraph $\Gamma_{\alpha}$ of $\alpha$ 
containing $[\omega]$
in its interior (see \cite{Massart}), and $\mbox{Vect} \, F_{\omega}$ the 
underlying vector space of the affine subspace generated by $F_{\omega}$
in $H^1(M, \R)$. 
Beware that $\mbox{Vect} \, F_{\omega}$ is not, unless $F_{\omega}$ contains 
the origin,, the vector space generated by $F_{\omega}$.
Note that $F_{\omega} = \{[\omega]\}$
if $\alpha$ is strictly convex at $[\omega]$. 
The value of $\alpha$ at the 
null cohomology class is Ma\~n\'e's critical value $c(L)$.

In section 3 we relate the dimension of the faces of $\Gamma_{\alpha}$
to the topological complexity of the complement of $\Pom$ in $M$, as follows.
Let $C_{\omega}(\epsilon)$ be the set of integer homology classes which are
represented by a piecewise $C^1$ closed curve 
made with arcs  contained  in \Pom , except for a remainder of total length 
less than $\epsilon$. 
Let $C_{\omega}$ be the 
intersection of $C_{\omega}(\epsilon)$ over all $\epsilon > 0$,
and let $V_{\omega}$ be the vector space spanned in 
$H_1(M, \R)$ by $C_{\omega}$. 
Note that $V_{\omega}$ is an integer subspace of $H_1(M, \R)$, that
is, it has a basis of integer elements (images in $H_1(M, \R)$ of elements of
$H_1(M, \Z)$).

We denote
\begin{itemize}
\item
by $V_{\omega}^{\perp}$ the vector
space of cohomology classes of one-forms  of class $C^1$
that vanish on $V_{\omega}$ .
\item
by $G_{\omega}$ the vector
space of cohomology classes of one-forms  of class $C^1$
that vanish in $T_x M$ for  
 every $x \in  \Pom $ 
\item
 by $E_{\omega}$ the vector
space of cohomology classes of one-forms of class $C^1$ , 
the supports of which are disjoint from \Pom .
\end{itemize}

\begin{theorem} 
\label{dimface}
We have  
$ E_{\omega} \subset \mbox{Vect}\, F_{\omega} \subset G_{\om}
\subset V_{\omega}^{\perp}$. 
When $M$ is a closed, orientable surface  all inclusions  are
equalities and  
furthermore  $\mbox{Vect}\, F_{\omega}$ is an integer subset of $H^1(M,\R)$.
\end{theorem}
\begin{theorem} \label{semicontinuite}
When $M$ is a closed, orientable surface
the vector space $\mbox{Vect}\, F_{\omega}$ is lower semi-continuous 
with respect to the  Lagrangian.
\end{theorem}

Theorem  \ref{dimface}  means that when $M$ is a closed, orientable surface,
the dimension of the face $F_{\omega}$ equals the number of homologically 
independant closed curves disjoint from $\Pom$.

As a corollary we get differentiability results for $\beta$.
The idea here was given to the author by Albert Fathi.

Let $h$ be a homology class. A cohomology class 
$\omega$ is said to be a subderivative for $\beta$ at $h$  
if $<\om, h> = \beta(h) +\al (\om)$. 
The subderivatives for $\beta$ at $h$ form a face $F_h$ of $\Gamma_{\alpha}$. 
By proposition \ref{pomface} the Aubry (resp. Mather)
sets for all the cohomology classes in the interior of this face coincide. 
We call that Aubry set (resp. Mather set), 
the Aubry set (resp. Mather set) of $h$, and denote it $\Ph$ (resp. $\Mh$).  

Recall that the tangent cone to the epigraph of $\beta$ at $h$
is the smallest cone in $H_1 (M,\R) \times \R$ with vertex 
$(h, \beta (h))$ and containing the epigraph of $\beta$.
We say that the $\beta$-function is differentiable at $h$ in the direction
$d$ if the tangent cone to the epigraph of $\beta$ at $h$ contains the 
affine subspace $h + \R d$. 

Thus we say the $\beta$-function is differentiable in $k$ 
directions at a homology class $h$ if the tangent cone at $h$ to the epigraph
of $\beta$ splits as a metric product of $\R^k$ and another cone which 
contains no straight line (affine subspace of dimension one).

We say a homology class $h$ is $k$-irrational if $k$ is the dimension of the 
smallest subspace of $H_1(M, \R)$ generated by integer classes and containing
$h$. 
In particular $1$-irrational means ``on a line with rational slope''  and
$\dim H_1(M, \R)$-irrational means completely irrational. 
We call rational any homology class of the form $1/n h$ where $n$ is an integer
and $h$ is the image in $H_1(M,\R)$ of an integer homology class. 
The integrality of $\mbox{Vect}\, F_{\omega}$ has the following consequence :

\begin{corollary} \label{betadiff}
Let $M$ be a closed  orientable surface, 
and $L$ be a  Lagrangian on $M$.
At a $k$-irrational homology class $h$ the $\beta$-function of $L$ 
is differentiable in at least $k$ directions.
\end{corollary}
This was conjectured, and proved in the torus case, by V. Bangert.
A similar result was proved for twist maps of the annulus by J. Mather
in \cite{Mather1}.  See also \cite{Delgado}.

In particular when $M$ is a closed, orientable surface,
$\beta$ is  differentiable in every direction
at a completely irrational class.
Rademacher's theorem says a convex function is differentiable almost 
everywhere but does not provide an explicit set of differentiability points. 
In \cite{BIK} a $C^r$ metric is constructed on a 
torus of dimension $8r +8$, such that its stable norm is not differentiable
in all directions at some completely irrational class. 

On the other hand if $\beta$ is differentiable in one (resp. no)  direction at
some homology class $h$, then $h$ must be 1-irrational (resp. zero).
Also note that at every non-zero class $\beta$ is differentiable in the radial
direction.

In the next section we investigate generic properties of Lagrangian systems.
We say a property is true for a generic 
Lagrangian if, given a Lagrangian $L$,
there exists a residual (countable intersection of open and dense subsets)
subset $\mathcal O$ of $C^{\infty}(M)$ 
such that the property holds for $L+f,
\forall f \in \mathcal O$. 
Ma\~n\'e (\cite{Mane2, CDI}) proved that for a generic
Lagrangian, there exists a unique minimising measure and put forth 
in \cite{Mane2} the 

\begin{conjecture}[Ma\~n\'e] \label{manefaible}
For a generic Lagrangian $L$ on a closed manifold $M$
there exist a dense open set $U_0$ of $H^1(M, \R)$  such that 
$\forall \om \in U_0 $, $\Mom (L)$ consists of a single periodic orbit, 
or fixed point.
\end{conjecture}

As an application of Theorems \ref{dimface}, 
\ref{semicontinuite}, and the results
of \cite{Massart} we prove this conjecture to be true when $M$ is a
closed, orientable surface.

\section{Preliminary results}
Recall that by a theorem of Fathi (\cite{Fathi5}, p. 104) 
there exists a pair of conjugate weak KAM solutions $(u_+, u_-)$ such that
$u_+$ and $u_-$ coincide only on \Pom. 
The main result of this section is
\begin{proposition} \label{preliminaire}
For every $\ep >0$ there exists an integrable, 
non-negative function $G_{\ep}$ on $M$ such that 
$G_{\ep}\inv (0) = \Po $ and for every
absolutely continuous  arc  $\ga \co [0,t] \longrightarrow M$ 
we have 
\begin{equation}\label{inegalite}
\int_0^t (L+c(L))(\ga, \gadot) dt \geq u_+(\ga(t))- u_+(\ga(0))
+ \int_0^t G_{\ep}(\ga (t))dt - \ep.
\end{equation}
\end{proposition}
\proof 
Since  $M$ is compact and the functions $h_t$ are
equilipschitz on $M \times M$ (\cite{Mather3}, see also \cite{Fathi5}, p. 105),
by Ascoli's theorem, for every $\ep >0$ there exists $T >0$ such that 
\[
\forall x, y \in M,\;
t \geq T  \Rightarrow h_t(x, y) \geq h(x, y) -c(L)t -\ep. 
\]
Take $T(\ep)$ to be the infimum of such $T$'s.

Let $\gamma \co \R_+ \rightarrow M$ be a $C^1$ arc.
Take $\ep >0$. Let $\chi_{\ep}$ be $\ep/\max(1,T(\epsilon))$ times
the characteristic function of the closed set
$(u_- - u_+)^{-1}([2\ep, +\infty[)$. 
We prove,  for all positive $t$, 
\begin{equation}
\int_0^t (L+c(L))(\ga, \gadot)(s) ds 
\geq u_+(\ga(t)) -  u_+(\ga(0)) + \int_0^t \chi_{\ep}(\ga(s)) ds -\ep . 
\end{equation}
The proposition follows by taking $G_{\ep}$ to be 
the upper bound of the functions $\chi_{\delta}$ over all $\delta \leq \ep$.

Define a sequence in $\R_+$ by 
$t_0 = 0 $ and $t_{i+1} =$ 
\[
\max \{ t \geq t_i \co
t- t_i \geq T(\epsilon)  
 \mbox{ and Leb} ([t_i, t]\cap \ga ^{-1}(\supp(\chi_{\ep}))) 
\leq T(\epsilon)\}
\]
where $\mbox{Leb}$ denotes Lebesgue measure on $\R$. 
Observe that $\ga(t_i) \in \newline
\supp(\chi_{\ep})$,
that $t_{i+1} - t_i \geq T(\ep)$, and that for all $x$ between $t_i$ and 
$t_{i+1}$
\[
\int_{t_i}^{x}\chi_{\ep} (\ga )(s) ds \leq 
 \ep\frac{\mbox{Leb}([t_i, x] \cap 
\ga ^{-1}(\supp(\chi_{\ep})))}{\max(1,T(\epsilon))} \leq \ep.
\]

We have, taking $t_n$ to be the last $t_i$ before $t$, 
\[
 \int_0^t (L +c(L))(\ga, \gadot)(s) ds     = 
\]
\[ 
\sum_{t_{i+1} \leq t } 
\left( \int_{t_{i}}^{t_{i+1}}L (\ga, \gadot)(s) ds  
+ c(L)(t_{i+1}- t_{i}) \right) 
\]
\[
+\int_{t_{n}}^{t}L(\ga, \gadot)(s) ds + c(L)(t - t_{n}) 
\]

thus, since $u_{\pm}$ are weak KAM solutions, and by the definitions of $h_t$,
and  $T(\ep)$, 

\begin{eqnarray*}
\int_0^t L (\ga, \gadot)(s) ds  + c(L)t
& \geq  &  
\sum_{t_{i+1} \leq t} 
h_{t_{i+1}- t_{i}}(\ga(t_{i}),\ga(t_{i+1})) \\
& &
+ c(L)(t_{i+1}- t_{i})   + u_+(\ga(t)) - u_+(\ga(t_{n})) 
\end{eqnarray*}
\begin{eqnarray*}
\hspace{-1.5 cm} & \geq  & 
\sum_{t_{i+1} \leq t} \left(h(\ga(t_{i}),\ga(t_{i+1}))
-\ep \right) 
+ u_+(\ga(t)) - u_+(\ga(t_{n})) \\
& \geq & 
  \sum_{t_{i+1} \leq t} 
\left( u_-(\ga(t_{i+1})) - u_+(\ga(t_{i})) - \ep \right) 
+ u_+(\ga(t)) - u_+(\ga(t_{n})) \\
& \geq & 
\sum_{i = 0}^{n}  
\left( u_-(\ga(t_{i})) - u_+(\ga(t_{i})) - \ep \right) 
+ u_+(\ga(t)) - u_+(\ga(0)) \\
& \geq  &  u_+(\ga(t))- u_+(\ga(0))+ \ep \sharp \{i/\; t_{i} \leq t\} \\
& \geq &   
\sum_{i = 1}^{n}\int_{t_{i-1}}^{t_i} \chi_{\ep}(\ga)(s)) ds
+ u_+(\ga(t))- u_+(\ga(0)) \\
& \geq & u_+(\ga(t))- u_+(\ga(0))  + \int_0^t\chi_{\ep}(\ga)(s) ds \;  - \ep. 
\end{eqnarray*}
\qed

\subsection{proof of Proposition \ref{pomface}}
The next proposition enables us to speak of the Aubry set of a face of the 
epigraph of $\al$, and therefore, of the Aubry set of a homology class.
\begin{proposition} \label{pomface}
If a cohomology class $[\om_1]$ belongs to 
the maximal face $F_{\om}$ of $\Gamma_{\al}$ containing 
$[\om]$ in its interior,
then $\Pom \subset {\mathcal A}_{\om_{1}}$. 
In particular, if $[\om_1]$ belongs to
the interior of $F_{\om}$, then $\Pom = {\mathcal A}_{\om_{1}}$. 
Conversely,
if two cohomology classes $\om$ and $\om_1$ are such that 
$\tilde{\Pom } \cap \tilde{{\mathcal A}}_{\om_1} \neq \emptyset$, then 
$\al(\om)= \al(a\om + (1-a)\om_1)$ for all $a \in [0,1]$, i.e. 
$\Gamma_{\alpha}$
has a face containing $\om$ and $\om_1$.
\end{proposition}
\proof
We can find $\om_2 \in F_{\om}$ and $a \in \, ]0,1[$ such that 
$\om = a \om_1 + (1-a)\om_2$.
By \cite{Fathi3}, the following property characterises \Pom  : 
\[ 
\forall x \in \Pom,\; \exists \, t_n \rightarrow +\infty,\mbox{ and }
\ga_n \co [0, t_n] \rightarrow M,\mbox{ such that } \ga(0)=\ga(t_n)=x 
\]
\[
\mbox{ and } \int_0^{t_n} (L - \om)(\ga_n, \gadot_n)(s) ds  + \al(\om)t
\rightarrow 0.
\]
Now $\om = a \om_1 + (1-a)\om_2$, and $\al(\om)= a\al(\om_1)+(1-a)\al(\om_2)$
since $[\om_1],[\om_2] \in F_{\om}$.
Therefore
\[
a\left[\int_0^{t_n} (L - \om_1)(\ga_n, \gadot_n)(s) ds  + 
\al(\om_1)t \right] 
\]
\[
+(1-a)\left[\int_0^{t_n}(L - \om_2)(\ga_n, \gadot_n)(s) ds  + \al(\om_2)t
\right] \rightarrow 0.
\]
Observe that both summands on the left are non-negative, for if $u_-$ is a 
weak KAM solution for $L-\om_i$, $i=1,2$, we have
\begin{eqnarray*}
\int_0^{t_n} (L - \om_i)(\ga_n, \gadot_n)(s) ds  
+ \al(\om_1)t & \geq &
u_-(\ga(t_n)) - u_-(\ga(0)) \\
& = & u_-(x) -u_-(x)  =  0
\end{eqnarray*}
\[
\mbox{hence } \int_0^{t_n} (L - \om_i)(\ga_n, \gadot_n)(s) ds  + \al(\om_1)t
\rightarrow 0 \mbox{ when } n \rightarrow \infty.
\]
Conversely, let $(\ga, \gadot)$ be an orbit in 
$\tilde{\Pom} \cap \tilde{{\mathcal A}}_{\om_1}$. 
We have
\[ 
\int_0^{t} [L - \om  + \al(\om)](\ga, \gadot)(s) ds =
u(\ga(t))-u(\ga(0)),
\]
\[ 
\int_0^{t} [L - \om_1  + \al(\om_1)](\ga, \gadot)(s) ds = 
u_1(\ga(t))-u_1(\ga(0)),
\]
 where $u$ (resp. $u_1$) is a weak KAM solution for 
$L-\om$ (resp. $L-\om_1$). 
Therefore 
\[ 
\int_0^{t} [L - (a\om +(1-a)\om_1)  + \al(a\om +(1-a)\om_1)]
(\ga, \gadot)(s) ds = 
\]
\[
a[u(\ga(t))-u(\ga(0))] +(1-a)[u_1(\ga(t))-u_1(\ga(0))] 
\]
\[
+ t[ \al(a\om +(1-a)\om_1)-a\al(\om) -(1-a)\al(\om_1))] .
\]
The first two summands on the right are bounded below, hence for the sum to
be bounded below, we must have 
$\al(\om) = a\al(\om) +(1-a)\al(\om_1)= \al(a\om + (1-a)\om_1)$, since by
convexity of $\al$, $\al(a\om + (1-a)\om_1) \leq a\al(\om) +(1-a)\al(\om_1)$.
\qed

\section{Faces of the epigraph}

{\em Proof  of} $G_{\omega} \subset V_{\omega}^{\perp}$.

It amounts to showing that a one-form in $G_{\omega}$ vanishes on 
$V_{\omega}$. 
Let $\om \in G_{\omega}$ and 
let $h$ be represented as in the definition of $V_{\omega}$ for some 
$\epsilon > 0$. 
Call $S$ the part of the curve representing $h$ which 
consists in segments of \Pom  , and $R$ the remainder.
Now  $ < [\om], h> = \int_S \om +\int_R \om$ where 
the first summand is zero, and the second summand can be bounded by
$C\epsilon$, where $C$ depends on $L$ and $\om$ only. 
The conclusion 
follows since $ \epsilon$ is arbitrarily small.
\qed 

{\em Proof  of} $\mbox{Vect} \, F_{\omega} \subset G_{\omega}$.

Take $\om_1,\om_2 \in F_{\omega}$. By Proposition
\ref{pomface} the Aubry sets for $L-\om_1$ and $L-\om_2$
coincide with $\Pom $. The weak KAM solutions 
$(u_+, u_-)$ are differentiable at every point of $\Pu$ (see \cite{Fathi1})
with  derivative the Legendre transform of the (well defined)
tangent vector. 
This derivative is Lipschitz and furthermore 
(see \cite{Fathi5}, p. 92) we have 
\begin{equation} \label{Holder}
 |u_{\pm} (\phi (y)) - u_{\pm} (\phi (x)) - 
\frac{\partial L - \om}{\partial v}(\phi (x),\gadot (0)) \circ 
D_x \phi (y-x)| \leq K\|y-x\|^2 
\end{equation}
where $\phi$ is a local chart on $M$, $x$ and $y$ are two points in the 
inverse image of $\Pom$  by the chart, $\gadot (0)$ is the tangent 
vector to $\Pom$ at $\phi (x)$, and $K$ only depends on the chart.  
So Whitney's extension theorem (\cite{Federer}, theorem 3.1.14)
allows us to 
take $\tilde{u}_1,\tilde{u}_2$ two $C^1$
functions, the derivatives of which coincide with that of 
$u_+^1$ and $u_+^2$ respectively along  \Pom . 
Replace 
$ \om_2 $ by $\om_2 + d\,\tilde{u}_1 - d\,\tilde{u}_2$.  
This one-form 
coincides with $\om_1$ in the tangent space to every point of $\Pom $  hence
the cohomology class $[\om_1 - \om_2]$ belongs to $G_{\omega}$. 
\qed

{\em Proof  of}  $ E_{\omega} \subset \mbox{Vect}\, F_{\omega}$

Assume, replacing if necessary $L$ by $L-\om$, that $\om=0$.
We actually prove a slightly stronger statement.
Call $\mbox{$\tilde{{\mathcal T}}_0$}$  the 
intersection with $\Po$ of the
union of Hausdorff limits, when $\eta$ tends to $0$,
of supports of $L- \eta$-minimising measures, 
and call $\To$ its  projection to $M$.   
Let $\eta$ be supported away from $\To$.

For starters we prove that there exists $\delta >0$ such that
for all $L+\delta \eta$-minimising measure $\mu$, for all
$(x,v)$ in $\spt (\mu)$, we have $\delta |\eta_x (v)| \leq G_1 (x)$
where $G_1$ comes from Proposition \ref{preliminaire}.

Indeed, assume otherwise. 
Then there exists a sequence  $\delta_n \longrightarrow 0$, 
$L+\delta_n \eta$-minimising measures $\mu_n$, and points $(x_n,v_n)$ 
in $\spt (\mu_n)$ such that for all $n$, we have 
\begin{equation} \label{delta_n}
\delta_n |\eta_{x_n} (v_n )| > G_1 (x_n).
\end{equation}
The sequence $(x_n,v_n)$ is bounded in $TM$ because the measures
$\mu_n$ sit in the energy levels $\alpha (\delta_n [ \eta])$.
So we may assume $(x_n, v_n) \longrightarrow (x,v)$.
Then we have $G_1 (x_n) \longrightarrow 0$ so by construction of 
$G_1$, $G_1 (x) = 0$ and $x \in \Po$. Besides, $(x,v)$
belongs to a Hausdorff limit point of the sequence of compact sets
$\spt (\mu_n)$ so $x \in \To$. But then for $n$ large enough,
since $\eta$ is supported outside $\To$, we should have
$\eta_{x_n} (v_n ) =0$ which contradicts Equation \ref{delta_n}
since $G_1$ is non-negative.

Therefore we see that for every orbit $\ga$ in the support of an 
$L+\delta\eta$-minimising measure $\mu$, by Equation
\ref{inegalite} we have  
\[
\int_0^t (L + c(L) \pm \delta\eta ) (\ga, \gadot)(s) ds   
\geq 
\]
\[ u_+ (\ga(t)) - u_+ (\ga(0)) + 
\int_0^t ( G_1 \pm \delta\eta)(\ga, \gadot)(s)  ds - 1
\]
so , by averaging and letting $t$ go to infinity,
\[
-c(L \pm \delta\eta) \geq \int (G_1 \pm \delta\eta) d\mu \; - c(L)
\]
whence, since $G_1 \pm \delta\eta$ is non-negative on the support of
$\mu$,
\[  
\al(0)=c(L) \geq c(L\pm \delta\eta) = \al(\pm \delta\eta).
\]
By convexity of $\al$, we have 
$2 \al (0) \leq \al (\delta\eta) + \al(-\delta\eta)$ so
we get $\al(0)=\al(\pm\delta \eta)$ which implies that 
$\pm\delta[\eta] 
 $ belong to $F_0$.
\qed

\subsection{The two-dimensional case}
We prove that when $M$ is a closed surface,
$ V_{\omega}^{\perp} \subset E_{\omega}$, thus 
proving all inclusions to be equalities. Since $ V_{\omega}$ is an
integer subset of $H_1 (M, \R)$ this implies that 
$\mbox{Vect}\,F_{\omega}$ is an integer subset of $H^1 (M, \R)$.

To that end we prove that there exists a neighborhood $U$ of $\Pom$ such that
every closed curve contained in $U$ has its homology class in $V_{\omega}$.
First let us show how this implies the equality.
If a 1-form $\al$ vanishes on every element of $V_{\om}$,
then there exists a function $f$ 
defined on $U$ such that the restriction of $\al$ to $U$ is equal to 
$df$. 
Extend $f$ to $M$, now $\al -df$ is cohomologous to $\al$ and 
supported away from $U$.
 
Assume the surface has genus greater than one, the genus one case being treated
by Bangert in \cite{Bangert1}, and assume 
our reference metric $g$ has negative curvature. 
By \cite{Boyland-Gole} every minimising orbit stays
within finite distance, in the universal cover $\tilde{M}$ of $M$, of a 
$g$-geodesic. 
In particular one can define the ends of a minimiser in the 
boundary at infinity of $\tilde{M}$. 
Call $\lambda$ the geodesic lamination
obtained from $\Pom$ by replacing each orbit by the corresponding geodesic.

From \cite{CB}, we know that each boundary component of a 
connected component of the complementary set
of $\lambda$ in $M$ is  either a closed leaf of $\lambda$, or a finite
sequence of non-closed leaves $\delta_1, \ldots \delta_n$ such that 
$\delta_i$ and $\delta_{i+1}$ are asymptotic ($i$ being in $\Z/n\Z$).

Therefore each boundary component of a 
connected component of the complementary set
of $\Pom$ in $M$ is  either a closed orbit in $\Pom$, or a finite
sequence of non-closed orbits $\delta_1, \ldots \delta_n$ such that 
$\delta_i$ and $\delta_{i+1}$ are asymptotic ($i$ being in $\Z/n\Z$).

Hence for each boundary component $\delta$ of a 
connected component $R$ of the complementary set
of $\Pom$ in $M$ there exists a neighborhood $V$ of $\delta$ in $R$
such that every arc contained in $V$, with its end on $\delta$, is homotopic,
with fixed ends, to an arc consisting of portions of $\delta$ and a remainder 
of length arbitrarily small (or no remainder at all if $\delta$ is a closed 
leaf). 
Now we just need to take 
$U$ such that $U \cap M \setminus \Pom$ is contained
in the union over all boundary components of $R$, and over all 
connected component of the complementary set of $\Pom$ in $M$, of such 
neighborhoods.
\qed

{\em Proof  of Corollary \ref{betadiff}}.

Let $h$ be a $k$-irrational homology class. Then the set of subderivatives 
to $\beta$ at $h$ form a face $F_h$ of $\Gamma_{\alpha}$. Furthermore 
$\beta$ is differentiable in $(\dim H_1 (M,\R)-\dim F_h)$ directions.
Take $\omega$ in the interior of the face$F_h$. 
We have $F_h \subset F_{\om}$ so $\dim G_{\om} \geq \dim F_h$. 

Then for every 
$\omega' \in G_{\om}$ we have $<\omega', h> =0$. Note that
$\{h \in H_1 (M,\R)  \co <\omega', h> =0 \; \forall
\omega' \in G_{\om} \}$ is an integer subset of  $H_1 (M,\R)$, of dimension
$\dim H_1 (M,\R)-\dim G_{\om}$.

Since $h$ is $k$-irrational this implies 
$\dim H_1 (M,\R)-\dim G_{\om} \geq k$ whence  
$\dim H_1 (M,\R)-\dim F_h \geq k$
which proves Corollary \ref{betadiff}.
\qed

{\em Proof of Theorem \ref{semicontinuite}}

Assume a sequence of Lagrangians $L_n$ converges, 
in the $C^2$-topology, to a $C^2$ Lagrangian $L$.

Let $(u^+_n, u^-_n)$ be conjugate pair of weak KAM solutions for $L_n$.
By \cite{Fathi5}, p. 88 the functions $(u^+_n, u^-_n)$ are equi-Lipschitz.
By Ascoli's theorem we may assume that $(u^+_n, u^-_n)$ converges to a pair
$(u^+, u^-)$ of Lipschitz functions.  Furthermore $u^{\pm} \prec L +c(L)$.
Take $x \in M$ and $t \in \R_+$. For every $n \in \N$ there exists
a $C^1$ path $\gamma_n \co [0,t] \longrightarrow M$ such that 
$\gamma_n (0) = x$ (resp. $\gamma_n (t) = x$) 
and 
\[
u^{\pm}_n(\ga_n(t))-u^{\pm}_n(\ga_n(0)) =  
\int_0^t (L_n +c(L_n))(\gamma_n,\dot{\gamma}_n)ds .
\]
Take $v \in T_x M$ a
limit point of $\dot{\gamma}_n (0)$ (resp. $\dot{\gamma}_n (t)$).
Then the extremal trajectory $\gamma \co [0,t] \longrightarrow M$ of the 
Lagrangian, with $\gamma (0) = x$ and $\dot{\gamma} (0) = v$ 
(resp. $\gamma (t) = x$ and $\dot{\gamma} (t) = v$)is a uniform limit of
$\gamma_n$ and so 
\[
u^{\pm}(\ga(t))-u^{\pm}(\ga(0)) =  
\int_0^t (L +c(L))(\gamma,\dot{\gamma})ds .
\]
This shows that $(u^+, u^-)$ are  weak KAM solutions for $L$. 
Then for every neighborhood $U$
of $\{ x\in M \co u^+ (x) =  u^-(x)\}$ there exists an $N \in \N$ such that
$\forall n \geq N,\; \{ x\in M  \co u^+_n (x) =  u^-_n(x)\} \subset U$.
Hence there exists an $N \in \N$ such that 
$\forall n \geq N,\; E_0 (L) \subset E_0 (L_n)$.
\qed

\section{On Generic Lagrangians}
From now on we assume $M$ to be a closed orientable surface.
We  begin with a 
\begin{lemma} \label{Sdense}
Let $L$ be a  Lagrangian on a closed orientable surface. 
The set $S(L)$ of subderivatives to $\beta$ 
at 1-irrational homology classes is dense in $H^1(M,\R)$.
\end{lemma} 

\proof Assume there exists an open set $U$ in $H^1(M,\R)$ such that 
$U \cap S(L) = \emptyset$. We may assume $U$ to be convex. 
Then the set  
$V = \{h \in H_1(M,\R) \co \exists \om \in U, <\om,h>= \al(\om)+\beta (h)\}$ 
is also convex. 
Call $H$ the vector space $V$ generates in $H_1(M,\R)$.
Then, since $V$ does not contain any 1-irrational class, 
the codimension of $H$ is at least one. 
Now $U=\cup_{h \in V} F_h$
so there exists $h \in V$ such that $\dim F_h \geq 1$. Such an $h$ is
at most $(\dim H_1 (M, \R)-1)$-irrational by Corollary \ref{betadiff}. 
Take $\om$ in the 
interior of $F_h$ ; we have $\vect F_{\om} = E_{\om}$ so there exists a 
closed curve $\ga$, such that $\Pom$ is disjoint from $\ga$. 
Furthermore, by semi-continuity of  $\vect F_{\om} = E_{\om}$, 
there exists a convex neighborhood $U_1$
of $\om$ in $U$ such that for all $\om'$ in $U_1$, $\Po (L-\om')$ 
is disjoint from $\ga$. In particular $H$ is contained
in the integer subspace defined by the equation $\inter([\ga],.])=0$.

Now assume by induction we have proved 
that for some $2 \leq k \leq \dim H_1 (M,\R)-2$
there exist $\om_k$ in $U$, a convex neighborhood $U_k$
of $\om_k$ in $U$, and closed curves $\ga_1 := \ga, \ldots \ga_k$ 
such that for all $\om'$ in $U_k$, $\Po (L-\om')$ 
is disjoint from $\ga_1 := \ga, \ldots \ga_k$. 
Likewise define $V_k$ to be the set of homology classes at which
elements of $U_k$ are subderivatives, and $H_k$ to be the vector space
generated by $V_k$.
Then  $H_k$ is contained
in the integer subspace defined by the equations $\inter
([\ga_i],.])=0$ for $i=1, \ldots k$ and 
the codimension of  $H_k$ is at least $k$. 
Assume the codimension of $H_k$ is
exactly $k$ ; then as previously $H_k$ is an integer subspace. Any open 
(in the induced topology) subset of such a subspace contains  a 
1-irrational class, an impossibility.
So the codimension of $H_k$ is at least $k+1$. Then, as previously, 
there exists $h_k \in V_k$ such that $\dim F_{h_k} \geq k+1$. Such an $h_k$ is
at most $(\dim H_1 (M, \R)-k-1)$-irrational by Corollary \ref{betadiff}.
Take $\om_{k+1}$ in the interior of $F_{h_k}$ ; 
we have $\vect F_{\om_{k+1}} = E_{\om_{k+1}}$ so there exists a 
closed curve $\ga_{k+1}$ homologically independent from 
$\ga_1 := \ga, \ldots \ga_k$, 
such that $\Pom_{k+1}$ is disjoint from $\ga_1, \ldots \ga_{k+1}$. 
Furthermore, by semi-continuity of  $\vect F_{\om} = E_{\om}$, 
there exists a convex neighborhood $U_{k+1}$
of $\om_{k+1}$ in $U_k$ such that for all $\om'$ in $U_{k+1}$, $\Po (L-\om')$ 
is disjoint from $\ga_1, \ldots \ga_{k+1}$.    

By induction we prove that $U$ contains a 
$(\dim H_1 (M,\R)-k)$-irrational class, for all 
$k = 1, \ldots \dim H_1 (M, \R)-1$, a contradiction.
\qed

By \cite{Massart}, Proposition 5, any minimizing  measure with a rational 
homology class must be supported on a union
of  periodic orbits, or  fixed points. 

By  \cite{Mane2}, Theorem D, for a given  homology class $h$, 
there exists a residual subset ${\mathcal O}_h$ of $C^{\infty} (M)$ 
such that for all $\phi \in {\mathcal O}_h$ 
there exists a unique closed measure in ${\mathcal M}_h (L+\phi)$. 

Then for all $h$ with rational direction,
for all $\phi \in {\mathcal O}_h$ 
there exists a unique closed measure 
$\mu_{h, \phi}$ in ${\mathcal M}_h (L+\phi)$,
supported on  a union of  periodic orbits $\ga_{h, \phi}$.
Every such periodic orbit is minimising in its homology class.
Then by \cite{Mane2}, Theorem D, 
we may assume that $\ga_{h, \phi}$ consists of 
pairwise non homologous periodic orbits. For any given $K$ 
only a finite number of integer homology classes have their 
$L$-action $\leq K$ so then 
$\ga_{h, \phi}$ actually consists of a finite number of periodic orbits
$\ga_{h, \phi, i}$.
For each of those orbits there exists a closed one-form $\omega_i$
such that $\ga_{h, \phi, i}$ is the unique $L-\omega_i$-minimising measure
(cf. \cite{Massart}, Theorem 8). 
Then by \cite{CI}, Theorem D, we may assume $ \ga_{h, \phi}$ to be 
hyperbolic in its energy level.

Next we prove  that, 
for all $\phi \in {\mathcal O}_{h}$, there exists $\epsilon (h,\phi)>0$, 
such that for any $\lambda  \in ]1-\epsilon (h,\phi), 1+\epsilon (h,\phi) [$, 
there exists a unique closed measure $\mu_{\lambda, \phi}$ 
in ${\mathcal M}_{\lambda h} (L+\phi)$, 
supported on a  union of  periodic orbits 
$\ga_{\lambda, \phi}$,  homotopic to $ \ga_{h, \phi}$.

Indeed, fix $\phi \in {\mathcal O}_{h}$, and 
consider a sequence $\lambda_n$ of real numbers converging to
one. Let  $\mu_n$ be $\lambda_n h$-minimising measures. The sequence of
measures $\mu_n$ converge to an  $h$-minimising measure, and the only
possibility is that it is supported on $ \ga_{h, \phi}$. The latter
being hyperbolic, a topological conjugacy argument proves our claim.

The set of  1-irrational homology classes 
is a countable union of lines. Choose  a countable dense subset 
$h_i,\; i \in \N$. Call ${\mathcal O}$  the intersection over all 
$i \in \N $ of ${\mathcal O}_{h_i}$ ;
this is a countable intersection of residual sets, hence residual.
Now for all $\phi \in {\mathcal O}$, there exists an open and dense subset $U(\phi)$ of the subset of  1-irrational homology classes, such that for any $h \in U(\phi)$,  
there exists a unique closed measure $\mu_{h, \phi}$ 
in ${\mathcal M}_{h} (L+\phi)$, supported on a union
of  periodic orbits $\ga_{h, \phi}$.

If $M$ has genus $\geq 2$, by Theorems 7 and 8 of \cite{Massart}, 
every subderivative to $\beta$ at a 1-irrational  homology class is 
contained in a face of codimension one, whether on the boundary 
or in the interior. By Corollary \ref{betadiff}, if a cohomology class is 
contained in a face of codimension one (resp. zero), then it must be 
subderivative to $\beta$ at a 1-irrational (resp. zero) homology class.

The same is true if $M$ is a torus and $\phi \in {\mathcal O}$ ; for
in that case, in every 1-irrational homology class $h$, there exists a
unique minimising measure. Such a measure is supported on one
periodic orbit, hence $\beta$ is not differentiable at $h$ 
(\cite{Bangert1}). 

Hence when $\phi \in {\mathcal O}$ $S(L+\phi)$ 
equals the set of cohomology class  
contained in a face of codimension one or zero.

Now consider the set $S'(L+\phi)$ of cohomology classes contained in the interior of
a face of codimension one or zero, and subderivative to $\beta$ at a point of $U(\phi)$. By Theorem \ref{semicontinuite} 
$S'(L+\phi)$ is open in $H^1(M,\R)$ for any Lagrangian $L$. 
Besides, since the interior of any face is dense in that face, and the $h_i$ are dense in $H_1 (m, \R)$,
$S'(L)$ is dense in $S(L)$, hence in $H^1(M,\R)$. 
Note that for all $\omega \in S'(L+\phi)$,  
$\Mom$ consists of periodic orbits with the same homology class, 
or fixed points. Indeed if 
$\Mom$ contained two homologically distinct periodic orbits, then  
$V_{\omega}$ would contain their homology classes and its dimension 
would be at least two, so $\omega$ could not lie in the the interior of a face
of codimension one or zero.

In particular
for for all $\phi \in {\mathcal O}$, $\omega \in S'(L+\phi)$, 
$\Mom (L+\phi)$ consists of one periodic orbit 
or  fixed point. This proves Conjecture \ref{manefaible} for surfaces.

{\small

{\bf Acknowledgements} : the author thanks 
Albert Fathi for patiently explaining his theory, 
Victor Bangert for his kind
invitation to Freiburg where part of this work was done, 
John Mather for pointing out a mistake in an early version,
Renato Iturriaga and Hector Sanchez 
for  the invitation to  the workshop on Lagrangian
Systems at CIMAT, and David Th\'eret, as well as Patrick Bernard, 
for useful conversations.

\bigskip

\noindent
CIMAT, Guanajuato, Gto., Mexico, and\\ 
GTA, UMR 5030, CNRS, Universit\'e Montpellier II, France\\
e-mail : massart@cimat.mx
}

\begin{thebibliography}{99}



\bibitem[Ba94]{Bangert1} Bangert, Victor  
{\em Geodesic rays, Busemann functions and monotone twist maps} 
Calc. Var. Partial Differential Equations 2 (1994), no.
1, 49--63.

\bibitem[Ba95]{Bangert2}Bangert, Victor 
{\em Minimal foliations and laminations}, 
Proceedings of the International Congress of Mathematicians, 
Vol. 1, 2 (Z\"urich, 1994), 453--464, Birkh\"auser, Basel, 1995.

\bibitem[Ba99]{Bangert3} Bangert, Victor 
{\em Minimal measures and minimizing closed normal one-currents}
 GAFA 9 (1999), no. 3, 413--427.

\bibitem[BG99]{Boyland-Gole} Boyland, Philip ; Gol\'e, Christophe 
{\em Lagrangian systems on hyperbolic
manifolds}, Ergodic Theory Dynam. Systems 19 (1999), no. 5, 1157--1173. 

\bibitem[BIK97]{BIK} Burago, D. ; Ivanov, S. ; Kleiner, B. 
{\em On the structure of the stable norm of periodic metrics}
 Math. Res. Lett. 4 (1997), no. 6,
791--808

\bibitem[C95]{Carneiro} Carneiro, Mario Jorge
{\em minimizing measures of the action of autonomous Lagrangians}, 
Nonlinearity 8 (1995), no. 6, 1077--1085.

\bibitem[CB88]{CB} Casson, Andrew J. ; Bleiler, Steven A. 
{\em Automorphisms of surfaces after Nielsen and Thurston}, 
London Mathematical Society Student Texts, 9. Cambridge University Press, 
Cambridge-New York, 1988. 

\bibitem[CDI97]{CDI}
G. Contreras, J. Delgado, R. Iturriaga 
{\em Lagrangian flows : the dynamics of globally minimizing orbits-II}, 
Bol. Soc. Brasil. Mat. (N.S.) {\bf 28} (1997), 155-196.

\bibitem[CI99]{CI} G. Contreras, R. Iturriaga,
{\em Convex Hamiltonian without conjugate points}, 
Ergodic Theory Dynam. Systems 19 (1999), no. 4, 901--952.



\bibitem[D93]{Delgado} J. Delgado, {\em Vertices of the action function of a
Lagrangian system}, Ph. D. thesis, IMPA, 1993.

\bibitem[Fa97a]{Fathi1}
A. Fathi 
{\em Th\'eor\`eme KAM faible et th\'eorie de Mather sur les syst\`emes
lagrangiens}, 
C. R. Acad. Sci. Paris, S\'erie I {\bf 324} (1997), 1043-1046.

\bibitem[Fa97b]{Fathi2}
A. Fathi {\em Solutions KAM faible et barri\`eres de Peierls},
 C. R. Acad. Sci. Paris, S\'erie I {\bf 325} (1997), 649-652.

\bibitem[Fa98a]{Fathi3}
A. Fathi {\em Orbites h\'et\'eroclines et ensemble de Peierls},
 C. R. Acad. Sci. Paris S\'er. I Math. 326 (1998), no. 10, 1213--1216. 

\bibitem[Fa98b]{Fathi4} A. Fathi 
{\em Sur la convergence du semi-groupe de Lax-Oleinik},   
C. R. Acad. Sci. Paris S\'er. I Math. 327 (1998), no. 3, 267--270.
%
\bibitem[Fa00]{Fathi5}
A. Fathi {\em Weak KAM Theorem in Lagrangian Dynamics}, to appear, Cambridge University Press.
%
\bibitem[Fe69]{Federer}
Federer, Herbert 
{\em Geometric measure theory}, Die Grundlehren der mathematischen
Wissenschaften, Band 153 Springer-Verlag New York Inc., New York 1969
%
\bibitem[Mn96]{Mane2} 
Ma\~n\'e,  Ricardo 
{\em Generic properties and problems of minimizing measures of 
Lagrangian systems} 
Nonlinearity 9 (1996), no. 2,
273--310
%
\bibitem[Mn97]{Mane3}
Ma\~n\'e,  Ricardo 
{\em Lagrangian flows : the dynamics of globally minimizing orbits-I}, \\
Bol. Soc. Brasil. Mat. (N.S.) {\bf 28} (1997), 141-155.
%
\bibitem[Mt97]{Massart} D. Massart 
{\em Stable norms for surfaces : local structure
of the unit ball at rational directions}, 
GAFA {\bf 7} (1997), 996-1010.
%
\bibitem[Mr90]{Mather1} Mather, John N. 
{\em  Differentiability of the minimal average action as a function of the
rotation number}, 
Bol. Soc. Brasil. Mat. (N.S.) 21 (1990), no. 1, 59--70.
% 
\bibitem[Mr91]{Mather2} J. N. Mather 
{\em Action minimizing invariant measures for positive definite  
Lagrangian systems}, Math. Z. {\bf 207}, 169-207 (1991).
%
\bibitem[Mr93]{Mather3} J. N. Mather 
{\em Variational construction of minimizing orbits},
Ann. Inst. Fourier {\bf 43}, 1349-1386 (1993).






\end{thebibliography}
\end{document}